# ESTIMATION OF THE HURST PARAMETER FROM DISCRETE NOISY DATA

By Arnaud Gloter and Marc Hoffmann

*Université de Marne-la-Vallée*

We estimate the Hurst parameter $H$ of a fractional Brownian motion from discrete noisy data observed along a high frequency sampling scheme. The presence of systematic experimental noise makes recovery of $H$ more difficult since relevant information is mostly contained in the high frequencies of the signal.

We quantify the difficulty of the statistical problem in a min-max sense: we prove that the rate $n^{-1/(4H+2)}$ is optimal for estimating $H$ and propose rate optimal estimators based on adaptive estimation of quadratic functionals.

**1. Introduction.**

1.1. *Motivation.* Many processes of interest in physics, molecular biology, finance and traffic networks possess, or are suspected to possess, self-similar properties. In this context, recovering the so-called scaling exponents from experimental data is a challenging problem. The purpose of this paper is to investigate a new statistical method for estimating self-similarity based on adaptive estimation of quadratic functionals of the noisy data by wavelet thresholding.

We stay with dimension 1 and focus on the paradigmatic example of fractional Brownian motion.

1.2. *Statistical model.* Let $X$ be a one-dimensional process of the form

$$X_t = \sigma W_t^H,$$

where $W^H$ is a fractional Brownian motion with self-similar index (or Hurst parameter) $H \in (0,1)$ and scaling parameter $\sigma \in (0,+\infty)$. In particular, $X$

---

is centered Gaussian with covariance $\mathbb{E}[X_s X_t]$ proportional to $|t|^{2H} + |s|^{2H} - |t-s|^{2H}$; see more in Section 4.1 below.

In practice, it is unrealistic to assume that a sample path of $X$ can be observed (in which case the parameters $H$ and $\sigma$ would be identified). Instead, $X$ is rather observed at discrete times. The problem of estimating $H$ and $\sigma$ in this context has been given considerable attention (some references are Dahlhaus [5], Istas and Lang [15] and Ludeña [18]).

In this paper we take the next logical step: we assume that each observation is contaminated by noise, so that for $i = 0, \ldots, n$ we observe

$$(1) \qquad Y_i^n = X_{i\Delta} + a(X_{i\Delta})\xi_i^n,$$

where the $\xi_i^n$ are (centered) noise terms and $\Delta^{-1}$ is the sampling frequency. The function $x \rightsquigarrow a(x)$ is an unknown nuisance parameter.

Throughout, we assume that the experiment lives over a fixed time horizon $[0,T]$, so we have $T = n\Delta$. With no loss of generality we take $T = 1$, hence $\Delta = \Delta_n = n^{-1}$. Recovering the Hurst parameter $H$ from the data $(Y_i^n)$ is our objective.

1.3. *Results.* We show in Theorems 1 and 2 below that the rate

$$v_n(H) = n^{-1/(4H+2)}$$

is optimal for estimating $H$. The accuracy $v_n(H)$ is slower by a polynomial order than the usual $n^{-1/2}$ obtained in the absence of noise. The difficulty lies in the fact that the information about $H$ is contained in the high frequencies of the signal $t \rightsquigarrow X_t$. Although the high frequency sampling rate $n$ usually allows one to recover $H$ at the classical rate $n^{-1/2}$ when $X$ is directly observed (by means of quadratic variations; see, e.g., [15]), the presence of the noise $\xi_i^n$ in this context significantly alters the nature of the problem.

## 2. Main results.

2.1. *Methodology.* The parameters $(H,\sigma)$ live in $\mathcal{D} \subset (0,1) \times (0,+\infty)$. The process $X$ and the noise variables $(\xi_i^n)$ are simultaneously defined on a common probability space endowed with a probability measure $\mathbb{P}_{H,\sigma}^n$ for each $n \geq 1$.

A rate $v_n \to 0$ is said to be achievable over $\mathcal{D}$ if there exists a (sequence of) estimator(s) $\widehat{H}_n$ such that the (sequence of) normalized error(s)

$$(2) \qquad v_n^{-1}(\widehat{H}_n - H)$$

is bounded in $\mathbb{P}_{H,\sigma}^n$-probability uniformly over $\mathcal{D}$. The rate $v_n$ is said to be a lower rate of convergence over $\mathcal{D}$ if there exists $c > 0$ such that

$$(3) \qquad \liminf_{n \to \infty} \inf_{\widehat{H}} \sup_{(H,\sigma) \in \mathcal{D}} \mathbb{P}_{H,\sigma}^n [v_n^{-1}|\widehat{H} - H| \geq c] > 0,$$

where the infimum is taken over all estimators $\widehat{H}$ that are random variables measurable with respect to the sigma-field generated by the data $(Y_i^n)$.

2.2. *The estimation strategy.* The fact that $X$ is a fractional Brownian motion enables one to predict that its energy levels

$$(4) \qquad Q_j := \sum_k d_{j,k}^2 := \sum_k \left( \int_{\mathbb{R}} X_s \psi_{j,k}(s)\, ds \right)^2$$

scale (as for the approximation symbol $\sim$, we do not yet specify it; see Proposition 1 below) with a ratio related to $H$,

$$(5) \qquad Q_{j+1} \sim 2^{-2H} Q_j,$$

up to an error term that vanishes as the frequency level $j$ increases. Here, $d_{j,k}$ is the random wavelet coefficient of the function $t \rightsquigarrow X_t$ relative to a certain wavelet basis $(\psi_{j,k}, j \geq 0, k \in \mathbb{Z})$. In Section 3.2 below we construct a procedure

$$(6) \qquad (Y_i^n) \rightsquigarrow (\widehat{d_{j,k,n}^2}, k = 0, \ldots, 2^j - 1, 0 \leq j \leq J_n)$$

that processes the data into estimates of the squared wavelet coefficients $d_{j,k}^2$ up to the maximal resolution level $J_n = [\frac{1}{2} \log_2(n)]$. We obtain a family of estimators for $H$ by setting

$$\widehat{H}_{j,n} := -\frac{1}{2} \log_2 \frac{\widehat{Q}_{j+1,n}}{\widehat{Q}_{j,n}}, \qquad j = 1, \ldots, J_n - 1,$$

with

$$\widehat{Q}_{j,n} = \sum_k \widehat{d_{j,k,n}^2}.$$

The ratio level $j$ between two estimated energy levels that contains maximal information about $H$ is chosen by means of a block thresholding rule; see below. The rule is inspired by the methodology introduced for the adaptive estimation of quadratic functionals (see, among others, Efromovich and Low [7], Gayraud and Tribouley [9] and the references therein).

2.3. *Statement of the results.* We consider for $(H, \sigma)$ regions of the form

$$(7) \qquad \mathcal{D} := [H_-, H_+] \times [\sigma_-, \sigma_+] \subset (\tfrac{1}{2}, 1) \times (0, +\infty).$$

ASSUMPTION A. (i) The function $x \rightsquigarrow a(x)$ is bounded and continuously differentiable with a bounded derivative.

(ii) The continuous time process $X$ is $\mathcal{F}^n$-adapted with respect to a filtration $\mathcal{F}^n = (\mathcal{F}_t^n, t \geq 0)$.

(iii) The noise term $\xi_i^n$ at time $i/n$ is $\mathcal{F}_{(i+1)/n}^n$-measurable. Moreover,

$$\mathbb{E}_{H,\sigma}^n[\xi_i^n|\mathcal{F}_{i/n}^n] = 0, \qquad \mathbb{E}_{H,\sigma}^n[(\xi_i^n)^2|\mathcal{F}_{i/n}^n] = 1,$$

and

$$\sup_{(H,\sigma)\in\mathcal{D}} \sup_{i,n} \mathbb{E}_{H,\sigma}^n[(\xi_i^n)^4] < +\infty.$$

THEOREM 1. *Grant Assumption* A. *The rate* $v_n(H) := n^{-1/(4H+2)}$ *is achievable for estimating* $H$ *over any region* $\mathcal{D}$ *of the form* (7). *Moreover, the estimator constructed in Section* 3 *and given by* (9)–(11) *below achieves the rate* $v_n(H)$.

This rate is indeed optimal as soon as the noise process enjoys some regularity:

ASSUMPTION B. (i) $\inf_x a(x) > 0$.

(ii) Conditional on $X$, the variables $\xi_i^n$ are independent, absolutely continuous with $\mathcal{C}^2$ densities $x \rightsquigarrow \exp(-v_{i,n}(x))$ vanishing at infinity (together with their derivatives) at a rate strictly faster than $1/x^2$ and

$$(8) \qquad \sup_{i,n} \mathbb{E}\left[\left(\frac{d}{dx}v_{i,n}(\xi_i^n)\right)^2 (1+|\xi_i^n|^2)\right] < +\infty.$$

Moreover, the functions $x \rightsquigarrow \frac{d^2}{dx^2}v_{i,n}(x)$ are Lipschitz continuous, with Lipschitz constants independent of $i,n$.

THEOREM 2. *Grant Assumptions* A *and* B. *For estimating* $H$, *the rate* $v_n(H) := n^{-1/(4H+2)}$ *is a lower rate of convergence over any region* $\mathcal{D}$ *of the form* (7) *with nonempty interior.*

We complete this section by giving an ancillary result about the estimation of the scaling parameter $\sigma$, although we are primarily interested in recovering $H$. The estimation of $\sigma$ has been addressed by Gloter and Jacod [12] for the case $H = 1/2$ and by Gloter and Hoffmann [10] in a slightly different model when $H \geq 1/2$ is known. Altogether, the rate $v_n(H)$ is proved to be optimal for estimating $\sigma$ when $H$ is *known*. Our next result shows that we lose a logarithmic factor when $H$ is *unknown*.

THEOREM 3. *Grant Assumptions* A *and* B. *For estimating* $\sigma$, *the rate* $n^{-1/(4H+2)}\log(n)$ *is a lower rate of convergence over any region of the form* (7).

### 2.4. *Discussion.*

2.4.1. *About the rate.* We see that the presence of noise dramatically alters the accuracy of estimation of the Hurst parameter: the optimal rate $v_n(H) = n^{-1/(4H+2)}$ inflates by a polynomial order as $H$ increases. In particular, the classical (parametric) rate $n^{-1/2}$ is obtained by formally letting $H$ tend to 0 (a case we do not have here).

2.4.2. *About Theorem* 1. The restriction $H_- > 1/2$ is linked to the discretization effect of the estimator. Assumption A can easily be fulfilled in the case of a noise process that is independent of the signal $X$. It is not minimal: more general noise processes could presumably be considered, and, more interestingly, scaling processes more general than fractional Brownian motion as well. To this end, it is required that the energy levels of $X$ satisfy Proposition 1 and that the empirical energy levels satisfy Proposition 2 in Section 4 below. We do not pursue that here. See also Lang and Roueff [17].

2.4.3. *About Theorem* 2. The lower bound is local, in the sense that $\mathcal{D}$ can be taken arbitrarily small in the class specified by (7). Observe that since the rate $v_n(H)$ depends on the parameter value, the min-max lower bounds (3) are only meaningful for parameter sets $\mathcal{D}$ that are concentrated around some given value of $H$.

Assumption B(ii) is not minimal: it is satisfied, in particular, when the $\xi_i^n$ are i.i.d. centered Gaussian. More generally, any noise process would yield the same lower bound as soon as Proposition 4 is satisfied (see Section 6.1).

2.4.4. *The stationary case.* Golubev [13] remarked that in the particular case of i.i.d. Gaussian noise independent of $W^H$, a direct spectral approach is simpler. Indeed, the observation generated by the $Y_i^n - Y_{i-1}^n$ becomes stationary Gaussian, and a classical Whittle estimator will do (Whittle [25] or Dahlhaus [5]). In particular, although some extra care has to be taken about the approximation in $n$, such an approach would certainly prove simpler in that specific context for obtaining the lower bound.

2.4.5. *Quadratic variation alternatives.* The estimator constructed in Section 3 can be linked to more traditional quadratic variation methods. Indeed, the fundamental energy levels $Q_j$ defined in (4) can be obtained from the quadratic variation of $X$ in the particular case of the Schauder basis (which does not have sufficiently many vanishing moments for our purpose). However, the choice of an optimal $j$ remains and we were not able to obtain the exact rate of convergence by this approach.

2.5. *Organization of the paper.* In Section 3 we give the complete construction of an estimator $\widehat{H}_n$ that achieves the min-max rate $v_n(H)$. Section 4 explores the properties of the energy levels of $X$ (Proposition 1), as well as their empirical version (Proposition 2). Theorem 1 is proved in Section 5. Finally, Sections 6 and 7 are devoted to the lower bounds. It is noteworthy that the complex stochastic structure of the model due to the two sources or randomness ($W^H$ and the noise $\xi_i^n$) requires particular efforts for the lower bound. Our strategy is outlined in Section 6: it requires a "coupling" result proved in Section 7. The proof of supplementary technical results, too long to be detailed here, may be found in [11].

## 3. Construction of an estimator.

3.1. Pick a wavelet basis $(\psi_{j,k}, j \geq 0, k \in \mathbb{Z})$ generated by a mother wavelet $\psi$ with two vanishing moments and compact support in $[0, S]$, where $S$ is some integer. The basis is fixed throughout Sections 3–5. Assuming we have estimators $\widehat{d^2_{j,k,n}}$ of the squared wavelet coefficients, recalling the definition (4) of the energy levels, we obtain a family of estimators for $H$ by setting

$$\widehat{H}_{j,n} := -\frac{1}{2}\log_2 \frac{\widehat{Q}_{j+1,n}}{\widehat{Q}_{j,n}}, \qquad j = \underline{J}, \ldots, J_n - 1,$$

with

$$\widehat{Q}_{j,n} = \sum_{k=0}^{2^{j-1}-1} \widehat{d^2_{j,k,n}},$$

where $J_n := [\frac{1}{2}\log_2(n)]$ is the maximum level of detail needed in our statistical procedure and $\underline{J} := [\log_2(S-1)] + 2$ is some (irrelevant) minimum level introduced to avoid border effects while computing wavelet coefficients corresponding to location on $[0, 1/2]$ from observations corresponding to $[0, 1]$. Following Gayraud and Tribouley [9] in the context of adaptive estimation of quadratic functionals, we let

(9) $$J_n^\star := \max\{j = \underline{J}, \ldots, J_n : \widehat{Q}_{j,n} \geq 2^j/n\}$$

(and in the case where the set above is empty, we let $J_n^\star = \underline{J}$ for definiteness). Eventually, our estimator of $H$ is

(10) $$\widehat{H}_{J_n^\star, n}.$$

The performance of $\widehat{H}_{J_n^\star, n}$ is related to scaling properties of $X$ and the accuracy of the procedure (6).

3.2. *Preliminary estimation of the $d_{j,k}^2$.* For simplicity and with no loss of generality, we assume from now on that $n$ has the form $n = 2^N$. Since $\psi$ has compact support in $[0, S]$, the wavelet coefficient $d_{j,k}$ is

$$d_{j,k} = \sigma \sum_{l=0}^{S2^{N-j}-1} \int_{k/2^j+l/2^N}^{k/2^j+(l+1)/2^N} \psi_{j,k}(t) W_t^H \, dt.$$

This suggests the approximation

$$\tilde{d}_{j,k,n} = \sum_{l=0}^{S2^{N-j}-1} \left( \int_{k/2^j+l/2^N}^{k/2^j+(l+1)/2^N} \psi_{j,k}(t) \, dt \right) Y_{k2^{N-j}+l}^n,$$

for $\underline{J} \leq j \leq J_n, 0 \leq k \leq 2^{j-1} - 1$. The difference $\tilde{d}_{j,k,n} - d_{j,k}$ splits into $b_{j,k,n} + e_{j,k,n}$, respectively a bias term and a centered noise term,

$$b_{j,k,n} = - \sum_{l=0}^{S2^{N-j}-1} \int_{k/2^j+l/2^N}^{k/2^j+(l+1)/2^N} \psi_{j,k}(t)(X_t - X_{k/2^j+l/2^N}) \, dt,$$

$$e_{j,k,n} = \sum_{l=0}^{S2^{N-j}-1} \left( \int_{k/2^j+l/2^N}^{k/2^j+(l+1)/2^N} \psi_{j,k}(t) \, dt \right) a(X_{k/2^j+l/2^N}) \xi_{k2^{N-j}+l}^n.$$

We denote by $v_{j,k,n}$ the variance of $e_{j,k,n}$, conditional on $\mathcal{F}_{k2^{-j}}^n$, which is equal to

$$v_{j,k,n} = \sum_{l=0}^{S2^{N-j}-1} \left( \int_{k/2^j+l/2^N}^{k/2^j+(l+1)/2^N} \psi_{j,k}(t) \, dt \right)^2 \mathbb{E}_{H,\sigma}^n[a(X_{k/2^j+l/2^N})^2 \mid \mathcal{F}_{k2^{-j}}^n].$$

The conditional expectations appearing in this expression are close to $a(X_{k/2^j})^2$ and thus may be estimated from the observations without the knowledge of $H, \sigma$. We define

$$\widehat{a^2}_{k/2^j,n} := 2^{-N/2} \sum_{l'=1}^{2^{N/2}} (Y_{k2^{N-j}+l'}^n)^2 - \left( 2^{-N/2} \sum_{l'=1}^{2^{N/2}} Y_{k2^{N-j}+l'}^n \right)^2$$

and we set

$$\overline{v}_{j,k,n} = \sum_{l=0}^{S2^{N-j}-1} \left( \int_{k/2^j+l/2^N}^{k/2^j+(l+1)/2^N} \psi_{j,k}(t) \, dt \right)^2 \widehat{a^2}_{k/2^j,n}.$$

Eventually, we set

(11) $$\widehat{d_{j,k,n}^2} := (\tilde{d}_{j,k,n})^2 - \overline{v}_{j,k,n}$$

and $\widehat{H}_{J_n^\star,n}$ is well defined. We remark that if the function $a$ is assumed known, one can considerably simplify the construction of the approximation $\widehat{a^2}_{k/2^j,n}$.

**4. The behavior of the energy levels.** We denote by $\mathbb{P}_{H,\sigma}$ the law of $X = \sigma W^H$, defined on an appropriate probability space. We recall the expression of the energy at level $j$,

$$Q_j = \sum_{k=0}^{2^{j-1}-1} d_{j,k}^2.$$

PROPOSITION 1. (i) *For all $\varepsilon > 0$, there exists $r_-(\varepsilon) \in (0,1)$ such that*

$$(12) \qquad \inf_{(H,\sigma) \in \mathcal{D}} \mathbb{P}_{H,\sigma} \left\{ \inf_{j \geq 1} 2^{2jH} Q_j \geq r_-(\varepsilon) \right\} \geq 1 - \varepsilon.$$

(ii) *The sequence*

$$(13) \qquad Z_j := 2^{j/2} \sup_{l \geq j} \left| \frac{Q_{l+1}}{Q_l} - 2^{-2H} \right|$$

*is bounded in $\mathbb{P}_{H,\sigma}$-probability, uniformly over $\mathcal{D}$, as $j \to +\infty$.*

PROPOSITION 2. *Let $j_n(H) := [\frac{1}{2H+1} \log_2(n)]$. Then $J_n \geq j_n(H)$ for all $H \in [H_-, H_+]$, and for any $L > 0$, the sequence*

$$n 2^{j_n(H)/2} \sup_{J_n \geq j \geq j_n(H) - L} 2^{-j} |\widehat{Q}_{j,n} - Q_j|$$

*is bounded in $\mathbb{P}_{H,\sigma}^n$-probability, uniformly over $\mathcal{D}$, as $n \to \infty$.*

We shall see below that Propositions 1 and 2 together imply Theorem 1.

4.1. *Fractional Brownian motion.* The fractional Brownian motion admits the harmonizable representation

$$W_t^H = \int_{\mathbb{R}} \frac{e^{it\xi} - 1}{(i\xi)^{H+1/2}} \mathrm{B}(d\xi),$$

where B is a complex Gaussian measure (Samorodnitsky and Taqqu [21]). Another representation using a standard Brownian motion $B$ on the real line is given by

$$W_t^H = \frac{1}{\Gamma(H + 1/2)} \int_{-\infty}^{\infty} [(t-s)_+^{H-1/2} - s_+^{H-1/2}] \, dB_s$$

($\Gamma$ is the Euler function). The process $W^H$ is $H$ self-similar and the covariance structure of $W^H$ is explicitly given by

$$\mathrm{Cov}(W_s^H, W_t^H) = \frac{\kappa(H)}{2} \{|t|^{2H} + |s|^{2H} - |t-s|^{2H}\},$$

where $\kappa(H) = \pi/H\Gamma(2H)\sin(\pi H)$. Recall that $d_{j,k} = \int_{\mathbb{R}} \psi_{j,k}(s) X_s \, ds$ denotes the random wavelet coefficients of $X$, given a wavelet $\psi$ with two vanishing moments. It can be seen, using the stationarity of the increments of $W^H$, that, for a fixed level $j$, the sequence $(d_{j,k})_{k \in \mathbb{Z}}$ is centered Gaussian and stationary with respect to the location parameter $k$. Moreover, the coefficients have the self-similarity property

$$(d_{j,k})_{k \in \mathbb{Z}} \stackrel{\text{law}}{=} 2^{-j(H+1/2)} (d_{0,k})_{k \in \mathbb{Z}};$$

see Delbeke and Abry [6], Veitch and Abry [23], Abry, Gonçalvès and Flandrin [1] and Veitch, Taqqu and Abry [24]. Moreover,

$$\operatorname{Var}(d_{j,k}) = \sigma^2 c(\psi) \kappa(H) 2^{-j(1+2H)},$$

where $c(\psi) = \frac{1}{2} \int \psi(s) \psi(t) \{|t|^{2H} + |s|^{2H} - |t-s|^{2H}\} \, ds \, dt$, and the covariance

$$\operatorname{Cov}(d_{j,k}, d_{j,k'}) = 2^{-j(2H+1)} \operatorname{Cov}(d_{0,k}, d_{0,k'})$$

decays polynomially as $k - k' \to \infty$ due to the two vanishing moments of $\psi$ and

$$|\operatorname{Cov}(d_{0,k}, d_{0,k'})| \leq c(1 + |k - k'|)^{2(H-2)},$$

for some $c$ which does not depend on $\sigma$ or $H$. See also Tewfik and Kim [22], Hirchoren and D'Attellis [14], Istas and Lang [15] and Gloter and Hoffmann [10].

PROPOSITION 3. *We have, for some constant $c > 0$,*

$$\sup_{(H,\sigma) \in \mathcal{D}} \mathbb{E}_{H,\sigma} \left[ \left( Q_j - 2^{-2jH} \frac{\sigma^2}{2} c(\psi) \kappa(H) \right)^2 \right] \leq c 2^{-j(1+4H)}.$$

PROOF. Remark that, by stationarity,

$$Q_j - 2^{-2jH} \frac{\sigma^2}{2} c(\psi) \kappa(H) = \sum_{k=0}^{2^{j-1}-1} (d_{j,k}^2 - \mathbb{E}_{H,\sigma}[d_{j,k}^2]).$$

Then the variance of the sum above is evaluated using the decorrelation property of the wavelet coefficients (similar computations can be found in Istas and Lang [15] and Gloter and Hoffmann [10]). □

4.2. *Proof of Proposition* 1. By Proposition 3, we derive in the same way as in Lemma II.4 of Ciesielski, Kerkyacharian and Roynette [4] that, for all $\varepsilon > 0$,

$$\sum_{j \geq 0} \sup_{(H,\sigma) \in \mathcal{D}} \mathbb{P}_{H,\sigma} \left[ 2^{2jH} Q_j \notin \left[ \frac{\sigma^2}{2} c(\psi) \kappa(H) - \varepsilon, \frac{\sigma^2}{2} c(\psi) \kappa(H) + \varepsilon \right] \right] < \infty,$$

from which (i) easily follows. By (i), the probability that $|Z_j|$ is greater than a constant $M$ is less than

(14) $$\varepsilon + \mathbb{P}_{H,\sigma}\left[\sup_{l \geq j}|Q_{l+1} - 2^{-2H}Q_l|2^{2lH} \geq M2^{-j/2}r_{-}(\varepsilon)\right].$$

By self-similarity, $\mathbb{E}_{H,\sigma}\{Q_{l+1} - 2^{-2H}Q_l\} = 0$. By Markov's inequality, (14) is less than

$$\varepsilon + [M^2 r_{-}^2(\varepsilon)]^{-1} \sum_{l \geq j} \text{Var}_{H,\sigma}(Q_{l+1} - 2^{-2H}Q_l)2^{4lH}2^j.$$

By Proposition 3, the sum above can be made arbitrarily small for large enough $M$, which proves (ii).

4.3. *Proof of Proposition* 2. We first claim that the following estimate holds:

(15) $$\sup_{J_n \geq j \geq j_n(H)-L} \sup_{(H,\sigma) \in \mathcal{D}} 2^{-j/2}\mathbb{E}^n_{H,\sigma}[|\widehat{Q}_{j,n} - Q_j|] \leq cn^{-1}.$$

Proposition 2 readily follows. To prove (15), we first split $\widehat{Q}_{j,n} - Q_j$ into $\sum_{u=1}^6 r_{j,n}^{(u)}$, with

$$r_{j,n}^{(1)} = \sum_k b_{j,k,n}^2, \qquad r_{j,n}^{(2)} = \sum_k (e_{j,k,n}^2 - v_{j,k,n}),$$

$$r_{j,n}^{(3)} = \sum_k (v_{j,k,n} - \overline{v}_{j,k,n}), \qquad r_{j,n}^{(4)} = 2\sum_k b_{j,k,n}d_{j,k},$$

$$r_{j,n}^{(5)} = 2\sum_k e_{j,k,n}d_{j,k}, \qquad r_{j,n}^{(6)} = 2\sum_k b_{j,k,n}e_{j,k,n}.$$

Using the result that $\mathbb{E}_{H,\sigma}[(X_t - X_s)^2] \leq c(H)\sigma^2|t-s|^{2H}$, it is readily seen that $\mathbb{E}^n_{H,\sigma}[(b_{j,k,n})^2]$ is less than a constant times $2^{-j}n^{-2H}$. Summing over $k$ shows that the term $r_{j,n}^{(1)}$ is negligible since $H > 1/2$.

Using the fact that $e_{j,k,n}^2 - v_{j,k,n}$ are uncorrelated for $|k - k'| \geq S$, we deduce that $\mathbb{E}^n_{H,\sigma}[(r_{j,n}^{(2)})^2]$ is bounded by a constant times $\sum_{k=0}^{2^{j-1}-1}\{\mathbb{E}^n_{H,\sigma}[e_{j,k,n}^4] + \mathbb{E}^n_{H,\sigma}[v_{j,k,n}^2]\}$. Then using the martingale increments structure of the sequence $a(X_{k2^{-j}+l2^{-N}})\xi^n_{k2^{-j}+l2^{-N}}$ for $l = 0, \ldots, S2^{N-j}$ (recall that $n = 2^N$), we may apply the Burkholder–Davis inequality. This gives, by Assumption A, $\mathbb{E}^n_{H,\sigma}[e_{j,k,n}^4] \leq cn^{-2}$. Then since $x \rightsquigarrow a(x)$ is bounded and, thus, $v_{j,k,n} \leq cn^{-1}$, we obtain that $\mathbb{E}^n_{H,\sigma}[(r_{j,n}^{(2)})^2]$ has the right order $2^j n^{-2}$.

Using conditional centering of $e_{j,k,n}$ with the fact that the variance of $d_{j,k}$ is less than $c2^{-j(2H+1)}$ and the condition $j \geq j_n(H) - L = [\frac{1}{2H+1}\log_2(n)] - L$, one easily checks that the terms $r_{j,n}^{(4)}$, $r_{j,n}^{(5)}$ and $r_{j,n}^{(6)}$ have negligible order.

We finally turn to the important term $r_{j,n}^{(3)}$, which encompasses the estimation of $a$. We claim that, for $0 \leq l \leq S2^{N-j}-1$, the following estimate holds:

$$(16) \quad \mathbb{E}_{H,\sigma}^n[|\widehat{a^2}_{k/2^j,n} - \mathbb{E}_{H,\sigma}^n[a(X_{k/2^j+l/2^N})^2 \mid \mathcal{F}_{k2^{-j}}^n]|] \leq cn^{-1/4}.$$

Summing over $l$ and $k$ yields the result for $r_{j,n}^{(3)}$ as soon as (16) is proved. Indeed, since $v_{j,k,n} - \overline{v}_{j,k,n}$ is equal to

$$\sum_{l=0}^{S2^{N-j}-1} \left( \int_{k/2^j+l/2^N}^{k/2^j+(l+1)/2^N} \psi_{j,k}(t) \, dt \right)^2 (\widehat{a^2}_{k/2^j,n} - \mathbb{E}_{H,\sigma}^n[a(X_{k/2^j+l/2^N})^2 \mid \mathcal{F}_{k2^{-j}}^n]),$$

we have that $\mathbb{E}_{H,\sigma}^n[|r_{j,n}^{(3)}|]$ is less than $c2^{j/2}n^{-1}2^{j/2}n^{-1/4}$. Therefore, under the restriction $j \leq J_n \leq [\frac{1}{2}\log_2(n)]$, (15) holds. It remains to prove (16).

We have $\widehat{a^2}_{k/2^j,n} - \mathbb{E}_{H,\sigma}^n[a(X_{k/2^j+l/2^N})^2 \mid \mathcal{F}_{k2^{-j}}^n] = t_{k,n}^{(1)} + t_{k,l,n}^{(2)} + t_{k,n}^{(3)}$, with

$$t_{k,n}^{(1)} = 2^{-N/2} \sum_{l'=1}^{2^{N/2}} X_{k/2^j+l'/2^N}^2 - \left( 2^{-N/2} \sum_{l'=1}^{2^{N/2}} Y_{k/2^{N-j}+l'}^n \right)^2,$$

$$t_{k,l,n}^{(2)} = 2^{-N/2} \sum_{l'=1}^{2^{N/2}} a(X_{k/2^j+l'/2^N})^2 (\xi_{k2^{j-N}+l'}^n)^2 - \mathbb{E}_{H,\sigma}^n[a(X_{k/2^j+l/2^N})^2 \mid \mathcal{F}_{k2^{-j}}^n],$$

$$t_{k,n}^{(3)} = 2^{-N/2+1} \sum_{l'=1}^{2^{N/2}} X_{k/2^j+l'/2^N} a(X_{k/2^j+l'/2^N}) \xi_{k2^{j-N}+l'}^n.$$

Since the $\xi_{k2^{j-N}+l'}^n$ are uncorrelated and centered, we readily have that the expectation of $|t_{k,n}^{(3)}|$ is of order $2^{-N/4} = n^{-1/4}$. For the term $t_{k,l,n}^{(2)}$, we use the preliminary decomposition

$$t_{k,l,n}^{(2)} = 2^{-N/2} \sum_{l'=1}^{2^{N/2}} a(X_{k/2^j+l'/2^N})^2 [(\xi_{k2^{j-N}+l'}^n)^2 - 1]$$

$$+ 2^{-N/2} \sum_{l'=1}^{2^{N/2}} (a(X_{k/2^j+l'/2^N})^2 - \mathbb{E}_{H,\sigma}^n[a(X_{k/2^j+l/2^N})^2 \mid \mathcal{F}_{k2^{-j}}^n]).$$

The expectation of the absolute value of the first term above is of order $n^{-1/4}$ since the summands $a(X_{k/2^j+l'/2^N})^2 [(\xi_{k2^{j-N}+l'}^n)^2 - 1]$ are martingale increments with second-order moments by Assumption A. Likewise, since $x \rightsquigarrow a(x)$ has a bounded derivative and

$$\mathbb{E}_{H,\sigma}[(X_{k/2^j+l'/2^N} - X_{k/2^j})^2] \leq c(H)\sigma^2(2^{-N/2})^{2H},$$

$$\mathbb{E}_{H,\sigma}\{(X_{k/2^j+l/2^N} - X_{k/2^j})^2\} \leq c(H)\sigma^2(2^{-j/2})^{2H},$$

the second term in the expression of $t^{(2)}_{k,l,n}$ has absolute expected value less than a constant times $(2^{-j/2})^H \leq 2^{j_n(H)H/2} = n^{-H/(1+2H)}$, and thus has the right order since $H \geq 1/2$.

Finally, we further need to split $t^{(1)}_{k,n}$ into

$$2^{-N/2} \sum_{l'=1}^{2^{N/2}} X^2_{k/2^j+l'/2^N} - \left(2^{-N/2} \sum_{l'=1}^{2^{N/2}} X_{k/2^j+l'/2^N}\right)^2$$

$$- \left(2^{-N/2} \sum_{l'=1}^{2^{N/2}} a(X_{k/2^j+l'/2^N})\xi^n_{k2^{N-j}+l'}\right)^2$$

$$- 2\left(2^{-N/2} \sum_{l'=1}^{2^{N/2}} X_{k/2^j+l'/2^N}\right)\left(2^{-N/2} \sum_{l'=1}^{2^{N/2}} a(X_{k/2^j+l'/2^N})\xi^n_{k2^{N-j}+l'}\right).$$

The first term and second term are easily seen to be of the right order, respectively, by the smoothness property of $X$ and the fact that the variables $\xi^n_i$ are uncorrelated. The third term is seen to have the right order after observing that one can replace the first sum $2^{-N/2}\sum_{l'=1}^{2^{N/2}} X_{k/2^j+l'/2^N}$ by $X_{k/2^j}$ up to a negligible error and then use the conditional zero correlation of the $\xi^n_i$ again. Thus, (16) is proved; hence, (15) follows. The proof of Proposition 2 is complete.

**5. Proof of Theorem 1.** First we need the following result that states the level $J^\star_n$, based on the data, is with large probability greater than some level based on the knowledge of $H$.

5.1. *A fundamental lemma.* For $\varepsilon > 0$, define

$$(17) \qquad J^-_n(\varepsilon) := \max\left\{j \geq 1; r^-(\varepsilon)2^{-2jH} \geq \frac{2^j}{n}\right\}.$$

LEMMA 1. *For all $\varepsilon > 0$, there exists $L(\varepsilon) > 0$ such that*

$$\sup_{(H,\sigma)\in\mathcal{D}} \mathbb{P}^n_{H,\sigma}[J^\star_n < J^-_n(\varepsilon) - L(\varepsilon)] \leq \varepsilon + \varphi_n(\varepsilon),$$

*where $\varphi_n$ satisfies $\lim_{n\to\infty} \varphi_n(\varepsilon) = 0$.*

PROOF. Let $L, \varepsilon > 0$. By definition of $J^-_n(\varepsilon)$,

$$\tfrac{1}{2}r_-(\varepsilon)^{1/(1+2H)}n^{1/(1+2H)} \leq 2^{J^-_n(\varepsilon)} \leq r_-(\varepsilon)^{1/(1+2H)}n^{1/(1+2H)};$$

hence, for large enough $n$, we have $\underline{J} \leq J_n^-(\varepsilon) - L \leq J_n$. Thus, by (9), $\mathbb{P}_{H,\sigma}^n[J_n^\star \geq J_n^-(\varepsilon) - L]$ is greater than

$$\mathbb{P}_{H,\sigma}^n[\widehat{Q}_{J_n^-(\varepsilon)-L,n} \geq 2^{J_n^-(\varepsilon)-L}n^{-1}],$$

which we rewrite as

$$\mathbb{P}_{H,\sigma}^n\{\widehat{Q}_{J_n^-(\varepsilon)-L,n} - Q_{J_n^-(\varepsilon)-L} \geq 2^{J_n^-(\varepsilon)-L}n^{-1} - Q_{J_n^-(\varepsilon)-L}\}$$

and which we bound from below by

$$\mathbb{P}_{H,\sigma}^n[\widehat{Q}_{J_n^-(\varepsilon)-L,n} - Q_{J_n^-(\varepsilon)-L} \geq 2^{J_n^-(\varepsilon)-L}n^{-1} - 2^{-2(J_n^-(\varepsilon)-L)H}r_-(\varepsilon)]$$
$$- \mathbb{P}_{H,\sigma}\left[\inf_{j \geq 1} 2^{2jH}Q_j < r_-(\varepsilon)\right].$$

Proposition 1(i) and the definition of $J_n^-(\varepsilon)$ yield that this last term is greater than

$$\mathbb{P}_{H,\sigma}^n[\widehat{Q}_{J_n^-(\varepsilon)-L,n} - Q_{J_n^-(\varepsilon)-L} \geq r_-(\varepsilon)^{1/(2H+1)}n^{-2H/(2H+1)}(2^{-L} - 2^{2LH})] - \varepsilon.$$

Then, if $L$ is such that $2^L - 2^{2LH} \leq -1$, an assumption we shall make from now on, Lemma 1 is proved, provided we show that

(18) $\quad \mathbb{P}_{H,\sigma}^n[|\widehat{Q}_{J_n^-(\varepsilon)-L,n} - Q_{J_n^-(\varepsilon)-L}| \geq r_-(\varepsilon)^{1/(2H+1)}n^{-2H/(2H+1)}]$

can be made arbitrarily small, uniformly in $(H,\sigma)$. Using again

$$2^{J_n^-(\varepsilon)} > \tfrac{1}{2}n^{1/(2H+1)}r_-(\varepsilon)^{1/(2H+1)},$$

we can pick $L' = L'(\varepsilon) > 0$ independent of $n$ such that

$$J_n^-(\varepsilon) - L \geq j_n(H) - L'(\varepsilon).$$

Therefore, (18) is less than

$$\mathbb{P}_{H,\sigma}^n\left[\sup_{J_n \geq j \geq j_n(H)-L'(\varepsilon)} |\widehat{Q}_{j,n} - Q_j| \geq r_-(\varepsilon)^{1/(2H+1)}n^{-2H/(2H+1)}\right],$$

which we rewrite as

$$\mathbb{P}_{H,\sigma}^n\left[n2^{j_n(H)/2} \sup_{J_n \geq j \geq j_n(H)-L'(\varepsilon)} 2^{-(j_n(H)-L'(\varepsilon))}|\widehat{Q}_{j,n} - Q_j| \geq v_H(\varepsilon,n)\right],$$

where

$$v_H(\varepsilon,n) := 2^{L'(\varepsilon)}r_-(\varepsilon)^{1/(2H+1)}n^{1/(4H+2)}$$

and where we use the fact that $2^{j_n(H)}$ is of order $n^{1/(2H+1)}$. We conclude by applying Proposition 2, using the fact that, for fixed $\varepsilon > 0$, $2^{L'(\varepsilon)}r_-(\varepsilon)^{1/(2H+1)}n^{1/(4H+2)} \to \infty$ as $n \to \infty$. The uniformity in $(H,\sigma)$ is straightforward. $\square$

5.2. *Proof of Theorem* 1, *completion.* Since $t \rightsquigarrow 2^{-2t}$ is invertible on $(0,1)$ with inverse uniformly Lipschitz on the compact sets of $(0,1)$, it suffices to prove Theorem 1 with $2^{-2H}$ in place of $H$ and $\widehat{Q}_{J_n^\star+1,n}/\widehat{Q}_{J_n^\star,n}$ in place of $\widehat{H}_{J_n^\star,n}$. First, we bound

$$\left|\frac{\widehat{Q}_{J_n^\star+1,n}}{\widehat{Q}_{J_n^\star,n}} - 2^{-2H}\right|$$

by a "bias" and a variance term, namely,

$$\left|\frac{Q_{J_n^\star+1}}{Q_{J_n^\star}} - 2^{-2H}\right| + \left|\frac{\widehat{Q}_{J_n^\star+1,n}}{\widehat{Q}_{J_n^\star,n}} - \frac{Q_{J_n^\star+1}}{Q_{J_n^\star}}\right| =: |B_n| + |V_n|,$$

say. Second, we prove Theorem 1 for $B_n$ and $V_n$ separately. We remark that the "bias" term $Q_{j+1}/Q_j - 2^{-2H}$ is deterministic, conditional on $X$, and decreases as the level $j$ increases, while the variance term $\widehat{Q}_{j+1,n}/\widehat{Q}_{j,n} - Q_{j+1}/Q_j$ increases. They both match at level $j = J_n^-(\varepsilon)$. In contrast to many "bias-variance" situations, the behavior of the variance term depends on the unknown regularity of the signal through the rate of decrease of the denominators $\widehat{Q}_{j,n}$ and $Q_j$. This explains the choice made in (9) to control the estimated level of energy $\widehat{Q}_{J_n^\star,n}$ from below.

5.2.1. *The bias term.* Let $M > 0$ and $\varepsilon > 0$. By Lemma 1, we have

$$\mathbb{P}^n_{H,\sigma}[n^{1/(4H+2)}|B_n| \geq M]$$
$$\leq \mathbb{P}^n_{H,\sigma}[n^{1/(4H+2)}|B_n| \geq M, J_n^\star \geq J_n^-(\varepsilon) - L(\varepsilon)] + \varepsilon + \varphi_n(\varepsilon)$$
$$\leq \mathbb{P}^n_{H,\sigma}[n^{1/(4H+2)} 2^{-J_n^-(\varepsilon)/2} 2^{L(\varepsilon)/2} |Z_{J_n^-(\varepsilon)-L(\varepsilon)}| \geq M] + \varepsilon + \varphi_n(\varepsilon)$$
$$\leq \mathbb{P}^n_{H,\sigma}[\sqrt{2} r_-(\varepsilon)^{-1/(4H+2)} 2^{L(\varepsilon)/2} |Z_{J_n^-(\varepsilon)-L(\varepsilon)}| \geq M] + \varepsilon + \varphi_n(\varepsilon),$$

where we have used for the last line the fact that, by (17),

$$2^{-J_n^-(\varepsilon)} \leq 2 r_-(\varepsilon)^{-1/(2H+1)} n^{-1/(2H+1)}.$$

We conclude by Proposition 1(ii) and by taking successively $\varepsilon$ sufficiently small, $M$ sufficiently large and $n$ sufficiently large.

5.2.2. *The variance term.* We split the variance term into $V_n = V_n^{(1)} + V_n^{(2)}$, where

$$V_n^{(1)} := \frac{\widehat{Q}_{J_n^\star+1,n} - Q_{J_n^\star+1}}{\widehat{Q}_{J_n^\star,n}} \quad \text{and} \quad V_n^{(2)} := \frac{Q_{J_n^\star+1}(Q_{J_n^\star} - \widehat{Q}_{J_n^\star,n})}{\widehat{Q}_{J_n^\star,n} Q_{J_n^\star}}.$$

Having Lemma 1 in mind, we bound, for any $M > 0$ and $L$ an integer, the probability $\mathbb{P}_{H,\sigma}^n[n^{1/(4H+2)}|V_n^{(1)}| \geq M]$ by

$$\mathbb{P}_{H,\sigma}^n[n^{1/(4H+2)}|V_n^{(1)}| \geq M, J_n^\star \geq J_n^-(\varepsilon) - L] + \mathbb{P}_{H,\sigma}^n[J_n^\star < J_n^-(\varepsilon) - L].$$

Fix $\varepsilon > 0$ and pick $L = L(\varepsilon)$ as in Lemma 1 so that the second probability $\mathbb{P}_{H,\sigma}^n[J_n^\star < J_n^-(\varepsilon) - L(\varepsilon)]$ is bounded by $\varepsilon + \varphi_n(\varepsilon)$. It remains now to deal with the first probability. As soon as $n$ is large enough, $J_n^-(\varepsilon) - L(\varepsilon) > \underline{J}$ and, thus, by definition of $J_n^\star$, the denominator of $V_n^{(1)}$ is bounded below by $2^{J_n^\star}/n$. This yields a new bound for the first probability,

$$\mathbb{P}_{H,\sigma}^n[n^{1/(4H+2)+1}2^{-J_n^\star}|\widehat{Q}_{J_n^\star+1,n} - Q_{J_n^\star+1}| \geq M, J_n^\star \geq J_n^-(\varepsilon) - L(\varepsilon)].$$

Recall that we defined $j_n(H) = [\frac{1}{2H+1}\log_2(n)]$ in Proposition 2 and by definition of $J_n^-(\varepsilon)$ we have

$$2^{J_n^-(\varepsilon)} > \tfrac{1}{2}n^{1/(2H+1)}r_-(\varepsilon)^{1/(2H+1)}.$$

Therefore, we can pick a positive $L' = L'(\varepsilon)$ independent of $n$ such that

$$J_n^-(\varepsilon) - L(\varepsilon) \geq j_n(H) - L'(\varepsilon),$$

and then we can bound the first probability by

$$\mathbb{P}_{H,\sigma}^n\left[n^{1/(4H+2)+1}\sup_{J_n \geq j \geq j_n(H) - L'(\varepsilon)} 2^{-j}|\widehat{Q}_{j,n} - Q_j| \geq M\right].$$

Next, using the fact that $n^{1/(4H+2)+1}$ is of order $n2^{j_n(H)/2}$ and Proposition 2, this term can be made arbitrarily small (uniformly in $n$) by taking $M$ large enough.

We now turn to the term $V_n^{(2)}$. Fix $\varepsilon > 0$ and $M > 0$. Recalling the definition of $Z_j$ in Proposition 1, we have

$$\mathbb{P}_{H,\sigma}^n[n^{1/(4H+2)}|V_n^{(2)}| \geq M]$$

$$\leq \mathbb{P}_{H,\sigma}^n\left\{n^{1/(4H+2)}\left|\frac{Q_{J_n^\star} - \widehat{Q}_{J_n^\star,n}}{\widehat{Q}_{J_n^\star,n}}\right|(2^{-2H} + Z_0) \geq M\right].$$

Now the tightness of the sequence $Z_j$ implies that, for some fixed constant $M'$, this probability is less than

$$\mathbb{P}_{H,\sigma}^n\left[n^{1/(4H+2)}\left|\frac{Q_{J_n^\star} - \widehat{Q}_{J_n^\star,n}}{\widehat{Q}_{J_n^\star,n}}\right| \geq \frac{M}{2^{-2H} + M'}\right] + \varepsilon.$$

Then the conclusion follows exactly as for $V_n^{(1)}$. The proof of Theorem 1 is complete.

**6. Proof of Theorems 2 and 3.** Consistently with Section 4, we denote by $\mathbb{P}_{H,\sigma}$ the probability measure on the Wiener space $\mathcal{C}_0$ of continuous functions on $[0,1]$ under which the canonical process $X$ has the law $\sigma W^H$. We write $\mathbb{P}_f^n$ for the law of the data, conditional on $X = f$.

6.1. *Preliminaries.* Define, for $\alpha \in (0,1)$ and $f \in \mathcal{C}_0$,

$$\|f\|_{\mathcal{H}^\alpha} := \|f\|_\infty + \sup_{0 \leq s < t \leq 1} \frac{|f(t) - f(s)|}{|t-s|^\alpha}, \tag{19}$$

with $\|f\|_\infty = \sup_t |f(t)|$.

The total variation of a signed measure $\mu$ is

$$\|\mu\|_{\mathrm{TV}} = \sup_{\|f\|_\infty \leq 1} \left| \int f \, d\mu \right|.$$

If $\mu$ and $\nu$ are two probability measures, the total variation of $\mu - \nu$ is maximal when $\mu$ and $\nu$ have disjoint support, in which case $\|\mu - \nu\|_{\mathrm{TV}} = 2$.

PROPOSITION 4. *Grant Assumptions* A *and* B. *We have, for some constant* $c > 0$,

$$\|\mathbb{P}_f^n - \mathbb{P}_g^n\|_{\mathrm{TV}} \leq c n^{1/2} \|f - g\|_\infty^{1/2}$$

*and*

$$1 - \tfrac{1}{2}\|\mathbb{P}_f^n - \mathbb{P}_g^n\|_{\mathrm{TV}} \geq R(cn\|f-g\|_2^2 + c\|f\|_{\mathcal{H}^{1/2}}^2 + c\|g\|_{\mathcal{H}^{1/2}}^2),$$

*where* $R$ *is some universal nonincreasing positive function and* $\|f\|_2 = (\int_0^1 f(s)^2)^{1/2}$.

PROOF. Let $D(\mu, \nu) := \int (\log \frac{d\mu}{d\nu}) d\mu \leq +\infty$ denote the Kullback–Leibler divergence between the probability measures $\mu$ and $\nu$. We recall the classical Pinsker inequality $\|\mu - \nu\|_{\mathrm{TV}} \leq \sqrt{2} D(\mu, \nu)^{1/2}$.

Using Assumption B(ii) and the representation (1), we deduce

$$\mathbb{E}_f^n\left[\log \frac{dP_f^n}{dP_g^n}(Y_0^n, \ldots, Y_n^n)\right]$$

$$= \sum_{i=0}^n \mathbb{E}_f^n\left[v_{i,n}(\xi_i^n + \Delta_{i,n}) - v_{i,n}(\xi_i^n) - \log\left(\frac{a(f_{i/n})}{a(g_{i/n})}\right)\right],$$

where $\Delta_{i,n} = \xi_i^n(\frac{a(f_{i/n})}{a(g_{i/n})} - 1) + \frac{f_{i/n} - g_{i/n}}{a(g_{i/n})}$. By a second-order Taylor expansion, this yields the expression for the Kullback–Leibler divergence,

$$D(dP_f^n, dP_g^n) = \sum_{i=0}^n \left\{ \mathbb{E}_f^n\left[\left(\frac{d}{dx} v_{i,n}\right)(\xi_i^n) \Delta_{i,n}\right] - \log\left(\frac{a(f_{i/n})}{a(g_{i/n})}\right) \right\}$$

(20)
$$+ \frac{1}{2} \sum_{i=0}^{n} \mathbb{E}_f^n \left[ \left( \frac{d^2}{dx^2} v_{i,n} \right)(\xi_i^n + \theta_{i,n}\Delta_{i,n})\Delta_{i,n}^2 \right],$$

for some (random) $\theta_{i,n} \in (0,1)$. Using the fact that $x \rightsquigarrow \exp(-v_{i,n}(x))$ vanishes at infinity, we have $\mathbb{E}_f^n[(\frac{d}{dx}v_{i,n})(\xi_i^n)] = 0$ and $\mathbb{E}_f^n[(\frac{d}{dx}v_{i,n})(\xi_i^n)\xi_i^n] = 1$, integrating by parts. It follows that the terms in the first sum of (20) are equal to $\frac{a(f_{i/n})}{a(g_{i/n})} - 1 - \log(\frac{a(f_{i/n})}{a(g_{i/n})})$. The assumptions on $x \rightsquigarrow a(x)$ yield that this quantity is less than some constant times $(f_{i/n} - g_{i/n})^2$.

For the second-order terms, using the uniform Lipschitz assumption on $x \rightsquigarrow \frac{d^2}{dx^2}v_{i,n}(x)$, together with the uniform bound for $\mathbb{E}_f^n[|\xi_i^n|^3]$, gives

$$\left| \mathbb{E}_f^n \left[ \left( \frac{d^2}{dx^2} v_{i,n} \right)(\xi_i^n + \theta_{i,n}\Delta_{i,n})\Delta_{i,n}^2 \right] \right|$$
$$\leq c|f_{i/n} - g_{i/n}|^3 + \left| \mathbb{E}_f^n \left[ \left( \frac{d^2}{dx^2} v_{i,n} \right)(\xi_i^n)\Delta_{i,n}^2 \right] \right|.$$

Again, we can bound $|\mathbb{E}_f^n[\frac{d^2}{dx^2}v_{i,n}(\xi_i^n)\Delta_{i,n}^2]|$ by a constant times $(f_{i/n} - g_{i/n})^2$, using the result that $\mathbb{E}_f^n[\frac{d^2}{dx^2}v_{i,n}(\xi_i^n)], \mathbb{E}_f^n[(\frac{d^2}{dx^2}v_{i,n}(\xi_i^n))\xi_i^n]$ and $\mathbb{E}_f^n[(\frac{d^2}{dx^2}v_{i,n}(\xi_i^n))(\xi_i^n)^2]$ are controlled by $\sup_{i,n} \mathbb{E}_f^n\{(\frac{d}{dx}v_{i,n}(\xi_i^n))^2(1 + |\xi_i^n|^2)\}$. Thus, the divergence between the conditional laws is bounded by

$$D(dP_f^n, dP_g^n) \leq c \sum_{i=0}^{n} |f_{i/n} - g_{i/n}|^2,$$

and the first part of the proposition follows from Pinsker's inequality. For the second part of the proposition, we use

$$\sum_{i=0}^{n} |f_{i/n} - g_{i/n}|^2 \leq 4n \int_0^1 (f(x) - g(x))^2 \, dx + 8n^{1-2\alpha}(\|f\|_{\mathcal{H}^\alpha}^2 + \|g\|_{\mathcal{H}^\alpha}^2),$$

valid for any $\alpha \in (0,1)$, together with the fact that for two measures $\mu, \nu$ the total variation $\|\mu - \nu\|_{\mathrm{TV}}$ remains bounded away from 2 when the divergences $D(\mu,\nu)$ and $D(\nu,\mu)$ are bounded away from $+\infty$. $\square$

The next result is the key to the lower bound. Its proof is delayed until Section 7. Let $(\sigma_0, H_0)$ be a point in the interior of $\mathcal{D}$. Set, for $I > 0$, $\varepsilon_n := I^{-1}n^{-1/(4H_0+2)}$ and

$$H_1 := H_0 + \varepsilon_n, \qquad \sigma_1 := \sigma_0 2^{j_0 \varepsilon_n},$$

where

$$j_0 = [\log_2(n^{1/(2H_0+1)})].$$

PROPOSITION 5. *For $I$ large enough, there exists a sequence of probability spaces $(\mathcal{X}^n, \mathfrak{X}^n, \mathbf{P}^n)$ on which can be defined two sequences of stochastic processes, $(\xi_t^{i,n})_{t \in [0,1]}$, $i = 0, 1$ such that:*

(i) *For $1/2 \leq \alpha < H_0$, the sequences $\|\xi^{0,n}\|_{\mathcal{H}^\alpha}$ and $\|\xi^{1,n}\|_{\mathcal{H}^\alpha}$ are tight under $\mathbf{P}^n$.*

(ii) *Define $P^{i,n} = \int_{\mathcal{X}^n} \mathbf{P}^n(d\omega) \mathbb{P}_{\xi^{i,n}(\omega)}^n$, and $\mathcal{Q}_{H,\sigma}^n = \int \mathbb{P}_{H,\sigma}(df) \mathbb{P}_f^n$, that is, the law of the data $(Y_i^n)$. Then*
$$\lim_{n \to \infty} \|P^{i,n} - \mathcal{Q}_{H,\sigma}^n\|_{\mathrm{TV}} = 0, \qquad i = 0, 1.$$

(iii) *There exists a measurable transformation $T^n: \mathcal{X}^n \mapsto \mathcal{X}^n$ such that the sequence $n\|\xi^{1,n}(\omega) - \xi^{0,n}(T^n(\omega))\|_2^2$ is tight under $\mathbf{P}^n$.*

(iv) *If $n$ is large enough, the probability measure $\mathbf{P}^n$ and its image measure $T^n \mathbf{P}^n$ are equivalent on $(\mathcal{X}^n, \mathfrak{X}^n)$. Moreover, for some $c^\star \in (0,2)$, we have*
$$\|\mathbf{P}^n - T^n \mathbf{P}^n\|_{\mathrm{TV}} \leq 2 - c^\star < 2,$$
*provided $n$ is taken large enough.*

REMARK. The processes $\xi^{0,n}$ and $\xi^{1,n}$ play the role of approximations for $\sigma_0 W^{H_0}$ and $\sigma_1 W^{H_1}$, respectively. Part (i) means that $\xi^{i,n}$ shares the same smoothness property as $W^{H_i}$, while (ii) implies that observing a noisy discrete sampling of $\sigma_i W^{H_i}$ ($i = 0, 1$) or of its approximation is statistically equivalent as $n \to \infty$. Of course, these points trivially hold in the case $\xi^{0,n} = \sigma_0 W^{H_0}$ and $\xi^{1,n} = \sigma_1 W^{H_1}$. However, a significant modification of this simple choice is needed in order to have the fundamental properties (iii) and (iv). These properties mean that one can transform pathwise, using $T^n$, the process $\xi^{0,n}$ into approximate realizations of $\xi^{1,n}$, while $T^n$ essentially does not transform $\mathbf{P}^n$ into a measure singular with it.

We next prove that Propositions 4 and 5 together imply Theorems 2 and 3.

6.2. *Proof of Theorems 2 and 3.* We prove Theorem 2 only. The proof of Theorem 3 is analogous since the choice of $H_i$ and $\sigma_i$ implies that $\sigma_1 - \sigma_0$ is of order
$$\frac{\sigma_0 \log(n)}{I(1 + 2H_0)} n^{-1/(2+4H_0)}.$$

Pick $n$ large enough so that $(\sigma_1, H_1) \in \mathcal{D}$. Pick an arbitrary estimator $\widehat{H}_n$. Let $M > 0$, with $M < 1/2I$ for further purposes. We have
$$\sup_{(H,\sigma) \in \mathcal{D}} \mathbb{P}_{H,\sigma}^n [n^{1/(4H+2)} |\widehat{H}_n - H| \geq M]$$

$$\geq \tfrac{1}{2}\mathbb{P}^n_{H_0,\sigma_0}[n^{1/(4H_0+2)}|\widehat{H}_n - H_0| \geq M]$$
$$+ \tfrac{1}{2}\mathbb{P}^n_{H_1,\sigma_1}[n^{1/(4H_1+2)}|\widehat{H}_n - H_1| \geq M]$$
$$\geq \tfrac{1}{2}P^{0,n}[n^{1/(4H_0+2)}|\widehat{H}_n - H_0| \geq M]$$
$$+ \tfrac{1}{2}P^{1,n}[n^{1/(4H_1+2)}|\widehat{H}_n - H_1| \geq M] + u_n,$$

where $u_n \to 0$ as $n \to \infty$ by (ii) of Proposition 5. By definition of $P^{i,n}$ and by taking $n$ large enough, it suffices to bound from below

$$(21) \qquad \tfrac{1}{2}\int_{\mathcal{X}^n}(\mathbb{P}^n_{\xi^0_n(\omega)}[A^0] + \mathbb{P}^n_{\xi^1_n(\omega)}[A^1])\mathbf{P}^n(d\omega),$$

where $A^i = \{n^{1/(2+4H_i)}|\hat{H}_n - H_i| \geq M\}$. By (iv) of Proposition 5, for $n$ large enough,

$$\int_{\mathcal{X}^n} \mathbb{P}^n_{\xi^{0,n}(\omega)}[A^0]\mathbf{P}^n(d\omega) = \int_{\mathcal{X}^n} \mathbb{P}^n_{\xi^{0,n}(\omega)}[A^0]\frac{d\mathbf{P}^n}{dT^n\mathbf{P}^n}(\omega)T^n\mathbf{P}^n(d\omega)$$
$$= \int_{\mathcal{X}^n} \mathbb{P}^n_{\xi^{0,n}(T^n(\omega))}[A^0]\frac{d\mathbf{P}^n}{dT^n\mathbf{P}^n}(T^n\omega)\mathbf{P}^n(d\omega).$$

Thus (21) is equal to half the quantity

$$\int_{\mathcal{X}^n}\left(\mathbb{P}^n_{\xi^{0,n}(T^n\omega)}[A^0]\frac{d\mathbf{P}^n}{dT^n\mathbf{P}^n}(T^n\omega) + \mathbb{P}^n_{\xi^{1,n}(\omega)}[A^1]\right)\mathbf{P}^n(d\omega)$$
$$\geq e^{-\lambda}\int_{\mathcal{X}^n}(\mathbb{P}^n_{\xi^{0,n}(T^n\omega)}[A^0] + \mathbb{P}^n_{\xi^{1,n}(\omega)}[A^1])1_{\frac{d\mathbf{P}^n}{dT^n\mathbf{P}^n}(T^n\omega)\geq e^{-\lambda}}\mathbf{P}^n(d\omega)$$
$$\geq e^{-\lambda}\int_{\mathcal{X}^n_r}(\mathbb{P}^n_{\xi^{0,n}(T^n\omega)}[A^0] + \mathbb{P}^n_{\xi^{1,n}(\omega)}[A^1])1_{\frac{d\mathbf{P}^n}{dT^n\mathbf{P}^n}(T^n\omega)\geq e^{-\lambda}}\mathbf{P}^n(d\omega),$$

for any $\lambda > 0$, and where $\mathcal{X}^n_r$ denotes the set of $\omega \in \mathcal{X}^n$ such that

$$n\|\xi^{0,n}(T^n\omega) - \xi^{1,n}(\omega)\|_2^2, \qquad \|\xi^{0,n}(T^n\omega)\|_{\mathcal{H}^\alpha} \quad \text{and} \quad \|\xi^{1,n}(\omega)\|_{\mathcal{H}^\alpha}$$

are bounded by $r > 0$. We will next need the two following technical lemmas.

LEMMA 2. *For any $r > 0$, there exists $c(r) > 0$ such that, on $\mathcal{X}^n_r$,*

$$\mathbb{P}^n_{\xi^{0,n}(T^n\omega)}[A^0] + \mathbb{P}^n_{\xi^{1,n}(\omega)}[A^1] \geq c(r) > 0.$$

LEMMA 3. *For large enough $n$, we have*

$$\mathbf{P}^n\left[\mathcal{X}^n_r \cap \frac{d\mathbf{P}^n}{dT^n\mathbf{P}^n}(T^n\cdot) \geq e^{-\lambda}\right] \geq \mathbf{P}^n[\mathcal{X}^n_r] - e^{-\lambda} - 1 + c^\star/2.$$

Applying successively Lemmas 2 and 3, we derive the lower bound

$$e^{-\lambda}c(r)(\mathbf{P}^n[\mathcal{X}_r^n] - e^{-\lambda} - 1 + c^\star/2).$$

Thus, Theorem 2 is proved as soon as we verify

(22) $$\lim_{r\to\infty}\liminf_{n\to\infty}\mathbf{P}^n[\mathcal{X}_r^n] = 1.$$

It suffices then to take $\lambda$ and $r$ large enough. By (i) and (iii) of Proposition 5, (22) only amounts to showing the tightness of $\|\xi^{0,n}(T^n\omega)\|_{\mathcal{H}^\alpha}$ under $\mathbf{P}^n$. For $L, L' > 0$, we have

$$\mathbf{P}^n[\|\xi^{0,n}(T^n(\omega))\|_{\mathcal{H}^\alpha} \geq L] = \int_{\mathcal{X}^n} 1_{\{\|\xi^{0,n}(\omega)\|_{\mathcal{H}^\alpha} \geq L\}} \frac{dT^n\mathbf{P}^n}{d\mathbf{P}^n}(\omega)\mathbf{P}^n(d\omega)$$

$$\leq L'\mathbf{P}^n[\|\xi^{0,n}(\omega)\|_{\mathcal{H}^\alpha} \geq L] + \mathbf{P}^n\left[\frac{dT^n\mathbf{P}^n}{d\mathbf{P}^n} \geq L'\right]$$

$$\leq L'\mathbf{P}^n[\|\xi^{0,n}(\omega)\|_{\mathcal{H}^\alpha} \geq L] + (L')^{-1}$$

by Chebyshev's inequality. The tightness of $\|\xi^{0,n}(T^n(\omega))\|_{\mathcal{H}^\alpha}$ then follows from the tightness of $\|\xi^{0,n}\|_{\mathcal{H}^\alpha}$. The proof of Theorem 2 is complete.

### 6.3. Proof of Lemmas 2 and 3.

#### 6.3.1. Proof of Lemma 2.

Since $H_0 < H_1$, it suffices to bound from below

$$\mathbb{P}^n_{\xi^{0,n}(T^n\omega)}[n^{1/(4H_0+2)}|\hat{H}_n - H_0| \geq M] + \mathbb{P}^n_{\xi^{1,n}(\omega)}[n^{1/(4H_0+2)}|\hat{H}_n - H_1| \geq M].$$

Let

$$d_{\text{test}}(\mu,\nu) := \sup_{0\leq f\leq 1}\left|\int f\,d\mu - \int f\,d\nu\right|$$

denote the test distance between the probability measures $\mu$ and $\nu$. The last term above is thus greater than

$$\mathbb{E}^n_{\xi^{1,n}(\omega)}[1_{n^{1/(4H_0+2)}|\hat{H}_n - H_0|\geq M} + 1_{n^{1/(4H_0+2)}|\hat{H}_n - H_1|\geq M}]$$
$$- d_{\text{test}}(\mathbb{P}^n_{\xi^{0,n}(T^n\omega)}, \mathbb{P}^n_{\xi^{1,n}(\omega)}).$$

Now since $M \leq 1/2I$ and by our choice for $H_0$ and $H_1$, one of the two events in the expectation above must occur with probability one. Using the fact that $d_{\text{test}}(\mu,\nu) = \frac{1}{2}\|\mu - \nu\|_{\text{TV}}$, the last term above is further bounded below by

$$1 - \tfrac{1}{2}\|\mathbb{P}^n_{\xi^{0,n}(T^n\omega)} - \mathbb{P}^n_{\xi^{1,n}(\omega)}\|_{\text{TV}}.$$

We conclude by Proposition 4 together with the fact that $\omega \in \mathcal{X}_r^n$.

6.3.2. *Proof of Lemma* 3. It suffices to bound from below

$$\mathbf{P}^n[\mathcal{X}_r^n] - \int_{\mathcal{X}^n} 1_{\frac{d\mathbf{P}^n}{dT^n\mathbf{P}^n}(T^n\omega) \leq e^{-\lambda}} \mathbf{P}^n(d\omega)$$

$$= \mathbf{P}^n[\mathcal{X}_r^n] - \int_{\mathcal{X}^n} 1_{\frac{dT^n\mathbf{P}^n}{d\mathbf{P}^n}(\omega) \geq e^{\lambda}} T^n\mathbf{P}^n(d\omega),$$

since $T^n\mathbf{P}^n$ and $\mathbf{P}^n$ are equivalent. We now replace the measure $T^n\mathbf{P}^n$ in the integral above by $\mathbf{P}^n$ with an error controlled by the test distance; the lower bounds become

$$\mathbf{P}^n[\mathcal{X}_r^n] - \mathbf{P}^n\left[\frac{dT^n\mathbf{P}^n}{d\mathbf{P}^n} \geq e^{\lambda}\right] - d_{\text{test}}(\mathbf{P}^n, T^n\mathbf{P}^n)$$

$$= \mathbf{P}^n[\mathcal{X}_r^n] - \mathbf{P}^n\left[\frac{dT^n\mathbf{P}^n}{d\mathbf{P}^n} \geq e^{\lambda}\right] - \frac{1}{2}\|\mathbf{P}^n - T^n\mathbf{P}^n\|_{\text{TV}}.$$

We conclude by the Chebyshev inequality and Proposition 5(iv).

**7. Proof of Proposition 5.** The proof of Proposition 5 relies on the construction of the fractional Brownian motion given by Meyer, Sellan and Taqqu [20]. In Section 7.1 we recall the main steps of the construction and how to apply it to our framework. In Section 7.2 we construct the sequence of spaces $(\mathcal{X}^n, \mathfrak{X}^n, \mathbf{P}^n)$. The proof of (i)–(iv) is delayed until Sections 7.3.1–7.3.4.

7.1. *A synthesis of fractional Brownian motion.* Consider a scaling function $\phi$ whose Fourier transform has compact support as in Meyer's book [19], with the corresponding wavelet function $\psi \in \mathcal{S}(\mathbb{R})$. In [20] the authors introduced, for $d \in \mathbb{R}$, the following differentials of order $d$ (via their Fourier transform):

$$\widehat{D^d\psi}(s) := (is)^d \hat{\psi}(s), \qquad \widehat{\phi^{d,\Delta}}(s) := \left(\frac{is}{1 - e^{is}}\right)^d \hat{\phi}(s),$$

where a determination of the argument on $\mathbb{C} \setminus \mathbb{R}_-$ with values in $(-\pi, \pi)$ is chosen. It is shown that the above formula is well defined and that $D^d\psi, \phi^{d,\Delta} \in \mathcal{S}(\mathbb{R})$. Define further, for $d = 1/2 - H \in (-1/2, 1/2)$,

$$\psi^H(t) := \int_{-\infty}^t D^d\psi(u)\,du = D^{d-1}\psi(t), \qquad \psi_{j,k}^H(t) := 2^{j/2}\psi^H(2^j t - k),$$

$$\Theta_k^H(t) := \int_0^t \phi^{d,\Delta}(u - k)\,du, \qquad \Theta_{j,k}^H(t) = 2^{j/2}\Theta_k^H(2^j t).$$

In their Theorem 2, Meyer, Sellan and Taqqu [20] prove the following almost sure representation of fractional Brownian motion (on an appropriate probability space and uniformly over compact sets of $\mathbb{R}$):

$$W_t^H = \sum_{k=-\infty}^{\infty} \Theta_k^H(t)\epsilon_k^H + \sum_{j=0}^{\infty}\sum_{k=-\infty}^{\infty} 2^{-j(H+1/2)}\{\psi_{j,k}^H(t) - \psi_{j,k}^H(0)\}\epsilon_{j,k},$$

where $\epsilon_k^H = \sum_{l=0}^{\infty} \gamma_l \epsilon_{k-l}'$ and $(1-r)^d = \sum_{k=0}^{\infty} \gamma_k r^k$ near $r=0$. The $\epsilon_k'$, $k \in \mathbb{Z}, \epsilon_{j,k}, j \geq 0, k \in \mathbb{Z}$ are i.i.d. $\mathcal{N}(0,1)$ random variables. Note that $\gamma_k = O(k^{-1+d})$, so the series above converges in quadratic mean and the time series obtained, $(\epsilon_k^H)_k$, has spectral density equal to $|2\sin(\frac{v}{2})|^{1-2H_0}$. The scaling

$$W_t^H \stackrel{\text{law}}{=} 2^{-j_0 H} W_{2^{j_0}t}^H$$

gives yet another representation for $W_t^H$,

$$\begin{aligned}(23)\quad &\sum_{k=-\infty}^{\infty} 2^{-j_0(H+1/2)}\Theta_{j_0,k}^H(t)\epsilon_k^H \\ &+ \sum_{j=j_0}^{\infty}\sum_{k=-\infty}^{\infty} 2^{-j(H+1/2)}\{\psi_{j,k}^H(t) - \psi_{j,k}^H(0)\}\epsilon_{j,k}.\end{aligned}$$

Comparing with other decompositions of fractional Brownian motion (e.g., Ciesielski, Kerkyacharian and Roynette [4] and Benassi, Jaffard and Roux [2]), a particular feature is that the random variables appearing in the high frequency terms

$$\sum_{j=j_0}^{\infty}\sum_{k=-\infty}^{\infty} 2^{-j(H+1/2)}\{\psi_{j,k}^H(t) - \psi_{j,k}^H(0)\}\epsilon_{j,k}$$

are independent and independent of the low frequency terms.

A drawback is that the basis used depends on $H$ and the functions appearing in the decomposition are not compactly supported. However, one can explore the properties of this basis. In [20], Meyer, Sellan and Taqqu show that the derivative of the initial wavelet function generates a multiresolution analysis and state the following results.

LEMMA 4 (Lemma 8 in [20]). *(1) There exist smooth $2\pi$-periodic functions $U_d$ and $V_d$ such that*

$$\widehat{\phi^{d,\Delta}}(s) = U_d(s/2)\widehat{\phi^{d,\Delta}}(s/2), \qquad \widehat{D^d\psi}(s) = V_d(s/2)\widehat{\phi^{d,\Delta}}(s/2).$$

*These "filters" and $U_d$ and $V_d$ vanish respectively in a neighborhood of $\pi$ and $0$.*

(2) Let $(c_k)_{k\in\mathbb{Z}} \in l^2(\mathbb{Z})$. Then the function $\sum_k c_k 2\phi^{d,\Delta}(2t-k)$ can be expressed with the basis $\phi^{d,\Delta}(t-k)$ and one level of detail,

$$\sum_k c_k 2\phi^{d,\Delta}(2t-k)$$
$$= \sum_k a_k \phi^{d,\Delta}(t-k) + \sum_k b_k D^d \psi(t-k), \tag{24}$$

where $(a_k)_{k\in\mathbb{Z}}$ and $(b_k)_{k\in\mathbb{Z}} \in l^2(\mathbb{Z})$. Moreover, $a$ and $b$ are given as follows: denoting by $A$, $B$ and $C$ the $2\pi$-periodic extensions of the discrete Fourier transforms of $a$, $b$ and $c$, we have

$$A(s) = -4^{-d}[V_d(s/2+\pi)C(s/2) - V_d(s/2)C(s/2+\pi)]e^{is/2}, \tag{25}$$

$$B(s) = -4^{-d}[-U_d(s/2+\pi)C(s/2) + U_d(s/2)C(s/2+\pi)]e^{is/2}. \tag{26}$$

From these properties we can show the following lemma, which will prove useful in controlling in $\mathcal{H}^\alpha$ norm the error made when we truncate the expansion. It also explores some properties of the basis when $H$ varies.

LEMMA 5. *Let $H \in (0,1)$. (i) If $u_k$ and $u_{j,k}$ are two sequences such that $|u_k| \leq c(1+|k|)^c$ and $|u_{j,k}| \leq c(1+j)^c(1+|k|)^c$, then, for any $\alpha \in [0,1)$ and $M \geq 0$, there exists $c(\alpha, M)$ such that, for all $j_0$,*

$$\sum_{j=j_0}^{\infty} \sum_{|k| \geq 2^{j+1}} \|u_{j,k}\psi^H_{j,k}\|_{\mathcal{H}_\alpha} \leq c(\alpha, M) 2^{-Mj_0},$$

$$\sum_{|k| \geq 2^{j_0+1}} \|u_k \Theta^H_{j_0,k}\|_{\mathcal{H}_\alpha} \leq c(\alpha, M) 2^{-Mj_0}.$$

(ii) *For all $M \geq 0$, there exists $c(M)$ such that, for all $\varepsilon > 0$ with $H+\varepsilon < 1$ and $t \in \mathbb{R}$,*

$$|\psi^{H+\varepsilon}(t) - \psi^H(t)| \leq c(M)\frac{\varepsilon}{(1+|t|)^M}. \tag{27}$$

(iii) *For all $\varepsilon > 0$ with $H+\varepsilon < 1$, we have, for all $k \in \mathbb{Z}$,*

$$\Theta^{H+\varepsilon}_k - \Theta^H_k$$
$$= \sum_{l\in\mathbb{Z}} a_l(\varepsilon)\Theta^H_{k+l} + \sum_{l\in\mathbb{Z}} b_l(\varepsilon)\{\psi^H_{0,k+l}(t) - \psi^H_{0,k+l}(0)\}, \tag{28}$$

*where the coefficients $a_l(\varepsilon)$ and $b_l(\varepsilon)$ are such that, for all $M$, there exists $c(M)$ such that, for all $\varepsilon$,*

$$\max\{|a_l(\varepsilon)|, |b_l(\varepsilon)|\} \leq \varepsilon c(M)(1+|l|)^{-M}. \tag{29}$$

*Moreover, the $2\pi$-periodic function $B_\varepsilon$ with Fourier coefficients $b_l(\varepsilon)$ vanishes in some neighborhood of zero independent of $\varepsilon$.*

The proof of Lemma 5 may be found in the Appendix of [11].

7.2. *The space* $(\mathcal{X}^n, \mathfrak{X}^n, \mathbf{P}^n)$. Let us recall that $H_1 = H_0 + \varepsilon_n$, where $\varepsilon_n = I^{-1} n^{-1/(2+4H_0)}$; $j_0 = [\log_2 n^{1/(1+2H_0)}]$ and $\sigma_1 = \sigma_0 2^{j_0 \varepsilon_n}$.

7.2.1. We take for $\mathcal{X}^n$ an infinite product of real lines, endowed with the product sigma field $\mathfrak{X}^n$,

$$\mathcal{X}^n := \left( \bigotimes_{k=-2^{j_0+1}}^{2^{j_0+1}} \mathbb{R} \right) \otimes \left( \bigotimes_{j=j_0}^{\infty} \bigotimes_{|k| \leq 2^{j+1}} \mathbb{R} \right) =: \mathcal{X}_\mathbf{e}^n \otimes \mathcal{X}_\mathbf{d}^n.$$

An element of $\mathcal{X}^n$ is denoted by $\omega = (\omega^\mathbf{e}, \omega^\mathbf{d})$ with $\omega^\mathbf{e} = (\omega_k^\mathbf{e})_{|k| \leq 2^{j_0+1}}$ and $\omega^\mathbf{d} = (\omega_\lambda^\mathbf{d})_{\lambda=(j,k); j \geq j_0, |k| \leq 2^{j+1}}$. The projections on the coordinates are denoted by $\epsilon_k(\omega) = \omega_k^\mathbf{e}$ for $|k| \leq 2^{j_0+1}$ and $\epsilon_{j,k}(\omega) = \omega_{j,k}^\mathbf{d}$ for $j \geq j_0, |k| \leq 2^{j+1}$.

On $\mathcal{X}^n$ we define the probability measure $\mathbf{P}^n := \mathbf{P}_\mathbf{e}^n \otimes \mathbf{P}_\mathbf{d}^n$, where $\mathbf{P}_\mathbf{e}^n$ is the unique probability on $\mathcal{X}_\mathbf{e}^n$ which makes the sequence $(\epsilon_k)$ a centered Gaussian stationary time series with spectral density $|2\sin(\frac{s}{2})|^{1-2H_0}$. The probability measure $\mathbf{P}_\mathbf{d}^n$ is the unique probability on $\mathcal{X}_\mathbf{d}^n$ that makes the sequence $(\epsilon_{j,k})$ i.i.d. $\mathcal{N}(0,1)$.

7.2.2. As suggested by Section 7.1, we define an approximation of $\sigma_0 W_0^H$ by keeping a finite number of coefficients at each scale,

$$\begin{aligned}
\xi^{0,n}(t) := & \sum_{|k| \leq 2^{j_0+1}} \sigma_0 2^{-j_0(H_0+1/2)} \Theta_{j_0,k}^{H_0}(t) \epsilon_k \\
& + \sum_{j \geq j_0} \sum_{|k| \leq 2^{j+1}} \sigma_0 2^{-j(H_0+1/2)} \{\psi_{j,k}^{H_0}(t) - \psi_{j,k}^{H_0}(0)\} \epsilon_{j,k}.
\end{aligned} \quad (30)$$

Denote by $\mathcal{T}^{n,1}$ a linear mapping from $\mathcal{X}_\mathbf{e}^n$ to itself such that, under the measure $\mathcal{T}^{n,1} \mathbf{P}_\mathbf{e}^n$, the coordinates $(\epsilon_k)$ form a centered Gaussian time series with spectral density $|2\sin(\frac{s}{2})|^{1-2H_1}$. Let

$$\epsilon_k'(\omega) := \epsilon_k(\mathcal{T}^{n,1}\omega). \quad (31)$$

We then define on the same space an approximation for $\sigma_1 W_1^H$. A natural choice would be to take again (30) with $(\sigma_1, H_1)$ and $\epsilon_k'$ instead of $(\sigma_0, H_0)$ and $\epsilon_k$. We proceed a little bit differently: we replace all the $\Theta_{j_0,k}^{H_1}$ by their truncated expansion on $\Theta_{j_0,k+l}^{H_0}$ and $\psi_{j_0,k+l}^{H_0}$ using relation (28). We then reorder the sums and finally drop the terms with index $k$ corresponding to the localization $k/2^j$ outside $[-2, 2]$. The reason is that we want to use the same basis as in $\xi^{0,n}$ for the low frequency terms.

This leads us to the following approximation for $\sigma_1 W^{H_1}$:

$$
\begin{aligned}
\xi^{1,n}(t) := & \sum_{|k|\leq 2^{j_0+1}} \sigma_1 2^{-j_0(H_1+1/2)} \Theta^{H_0}_{j_0,k}(t) \epsilon'_k \\
& + \sum_{|l|\leq 2^{j_0+1}} \sigma_1 2^{-j_0(H_1+1/2)} \Theta^{H_0}_{j_0,l}(t) \sum_{|k|\leq 2^{j_0+1}} a_{l-k}\epsilon'_k \\
& + \sum_{|l|\leq 2^{j_0+1}} \sigma_1 2^{-j_0(H_1+1/2)} \{\psi^{H_0}_{j_0,l}(t) - \psi^{H_0}_{j_0,l}(0)\} \sum_{|k|\leq 2^{j_0+1}} b_{l-k}\epsilon'_k \\
& + \sum_{j\geq j_0} \sum_{|k|\leq 2^{j+1}} \sigma_1 2^{-j(H_1+1/2)} \{\psi^{H_1}_{j,k}(t) - \psi^{H_1}_{j,k}(0)\} \epsilon_{j,k},
\end{aligned}
\tag{32}
$$

where the coefficients $a = a(\varepsilon)$ and $b = b(\varepsilon)$ are defined by (28) with $H = H_0$, $H + \varepsilon = H_1$.

7.2.3. The last step is the construction of the mapping $T^n$ from $(\mathcal{X}^n, \mathfrak{X}^n)$ to itself. Recalling (iii) of Proposition 5, we see that $T^n$ should transform outcomes of $\xi^{0,n}$ into approximate outcomes of $\xi^{1,n}$. Thus, we define the action of $T^n$ on the random space $(\mathcal{X}^n, \mathfrak{X}^n)$ by making the low frequency terms of $\xi^{0,n}(T^n\omega)$ exactly match the low frequency terms of $\xi^{1,n}(\omega)$.

We define $\mathcal{T}^{2,n}$ on $\mathcal{X}^n$ as the linear map such that

$$\epsilon_l(\mathcal{T}^{2,n}\omega) = \sum_{|k|\leq 2^{j_0+1}} a_{l-k}\epsilon_k(\omega) + \epsilon_l(\omega), \tag{33}$$

$$\epsilon_{j_0,l}(\mathcal{T}^{2,n}\omega) = \sum_{|k|\leq 2^{j_0+1}} b_{l-k}\epsilon_k(\omega) + \epsilon_{j_0,l}(\omega), \tag{34}$$

$$\epsilon_{j,l}(\mathcal{T}^{2,n}\omega) = \epsilon_{j,l}(\omega) \qquad \text{if } j > j_0. \tag{35}$$

We remark that the matrix of this linear map in the canonical basis of $\mathcal{X}^n$ is, of course, infinite, but $\mathcal{T}^{2,n}$ leaves invariant the finite-dimensional subspace $\mathcal{X}^n_{\mathbf{e}} \otimes (\otimes_{|k|\leq 2^{j_0+1}} \mathbb{R}) \otimes (0, 0, \ldots) \subset \mathcal{X}^n$ and is the identity on a supplementary space. On the finite-dimensional subspace its matrix is $\mathrm{Id} + K^n$, where $K^n$ is the square matrix of size $2[2^{j_0+2} + 1]$,

$$K^n = \begin{pmatrix} (a_{l-k})_{|l|,|k|\leq 2^{j_0+1}} & 0 \\ (b_{l-k})_{|l|,|k|\leq 2^{j_0+1}} & 0 \end{pmatrix}. \tag{36}$$

Finally, we set

$$T^n = \mathcal{T}^{n,2} \circ \mathcal{T}^{n,1}, \tag{37}$$

where we denote again by $\mathcal{T}^{n,1}$ the extension of $\mathcal{T}^{n,1}$ (previously defined only on $\mathcal{X}^n_{\mathbf{e}}$) to $\mathcal{X}^n$ such that it is the identity on $0_{\mathcal{X}^n_{\mathbf{e}}} \otimes \mathcal{X}^n_{\mathbf{d}}$.

As announced, the choice of $T^n$, with (30)–(35) and the fact that $\sigma_1 2^{-j_0 H_1} = \sigma_0 2^{j_0 \varepsilon_n} 2^{-j_0 H_1} = \sigma_0 2^{-j_0 H_0}$, yields

$$
\begin{aligned}
&\xi^{1,n}(\omega) - \xi^{0,n}(T^n(\omega)) \\
(38) \quad &= \sum_{j \geq j_0} \sum_{|k| \leq 2^{j+1}} \sigma_1 2^{-j(H_1+1/2)} \{\psi_{j,k}^{H_1}(t) - \psi_{j,k}^{H_1}(0)\} \epsilon_{j,k}(\omega) \\
&\quad - \sum_{j \geq j_0} \sum_{|k| \leq 2^{j+1}} \sigma_0 2^{-j(H_0+1/2)} \{\psi_{j,k}^{H_0}(t) - \psi_{j,k}^{H_0}(0)\} \epsilon_{j,k}(\omega).
\end{aligned}
$$

We now have completed the setup of $(\mathcal{X}^n, \mathfrak{X}^n, \mathbf{P}^n)$ and it now remains to prove that Proposition 5 holds. Let us stress that the choice of $j_0$ is for that matter crucial. Clearly, Proposition 5(iii) requires that $j_0$ be large enough. Meanwhile, Proposition 5(iv) requires that the number of components of $\mathcal{X}^n$ on which $T^n$ is different from the identity be as small as possible, which requires that $j_0$ be not too large. Since the proof is rather technical and quite long, we only sketch it here. A detailed proof may be found in [11].

### 7.3. Sketch of the proof of Proposition 5.

7.3.1. *Property* (i). We see that the representation (23) and our choice (30) only differ by the terms corresponding to locations $k/2^j \notin [-2, 2]$. With the help of Lemma 5(i), it can be deduced that on some probability space we have $\|\xi^{0,n} - \sigma_0 W^{H_0}\|_{\mathcal{H}_\alpha} \leq c(\omega) 2^{-M j_0}$, where $M$ is arbitrarily large and $c(\omega)$ is some random variable with finite moments coming from the randomness of the coefficient in the expansion (23). A similar bound may be obtained for $\|\xi^{1,n} - \sigma_1 W^{H_1}\|_{\mathcal{H}_\alpha}$. Then the property (i) of Proposition 5 follows from the almost sure smoothness property of the fractional Brownian motion.

7.3.2. *Property* (ii). Proposition 4 gives immediately an almost sure relation on the conditional laws: $\|\mathbb{P}^n_{\xi^{i,n}} - \mathbb{P}^n_{\sigma_i W^{H_i}}\|_{\mathrm{TV}} \leq c n^{1/2} \|\xi^{i,n} - \sigma_i W^{H_i}\|_\infty^{1/2}$ for $i = 0, 1$. Combining with the study of the difference $\xi^{i,n} - \sigma_i W^{H_i}$, this shows that this total variation distance is bounded by $c(\omega) n^{1/2} 2^{-M j_0/2}$. We are then able to deduce that the same bound holds for the unconditional laws

$$\|P^{i,n} - \mathcal{Q}^n_{H,\sigma}\|_{\mathrm{TV}} \leq c n^{1/2} 2^{-M j_0/2}.$$

Since $M$ is arbitrarily large, property (ii) of Proposition 5 follows and it is clear that this property is not crucial for the calibration of $j_0$.

7.3.3. *Property* (iii). We write (38) as

$$\xi^{1,n}(\omega) - \xi^{0,n}(T^n(\omega)) = q_1(t) - q_1(0) + q_2(t) - q_2(0),$$

where

$$q_1(t) := \sum_{j \geq j_0} \sum_{|k| \leq 2^{j+1}} \sigma_1 2^{-j(H_1+1/2)} \{\psi_{j,k}^{H_1}(t) - \psi_{j,k}^{H_0}(t)\} \epsilon_{j,k}(\omega),$$

$$q_2(t) := \sum_{j \geq j_0} \sum_{|k| \leq 2^{j+1}} (\sigma_1 2^{-j(H_1+1/2)} - \sigma_0 2^{-j(H_0+1/2)}) \psi_{j,k}^{H_0}(t) \epsilon_{j,k}(\omega).$$

But Lemma 5(ii) implies that the difference $\psi_{j,k}^{H_1} - \psi_{j,k}^{H_0}$ is a function with uniform norm bounded by $c2^{j/2}\varepsilon$ and well localized around $k/2^j$. This enables us to evaluate the sum with respect to $k$ in $q_1(t)$ and to deduce (for precise computations, see [11])

$$q_1(t)^2 \asymp \sum_{j \geq j_0} 2^{-2jH_1} \varepsilon^2 \asymp 2^{-2j_0 H_1} \varepsilon^2 \leq 2^{-2j_0 H_0} \varepsilon^2,$$

where $\asymp$ means equality in stochastic order. An analogous evaluation is obtained for $q_2(t)$, using $\sigma_1 2^{-jH_1} - \sigma_0 2^{-jH_0} = \sigma_0 2^{-jH_0}(2^{(j_0-j)\varepsilon} - 1)$.

Hence, property (iii) of Proposition 5 follows from $j_0 = [\frac{1}{2H_0+1} \log_2 n]$, which implies that $2^{-2j_0 H_0} \varepsilon^2$ is of order $n^{-1}$.

7.3.4. *Property* (iv). Let us focus only on the really delicate part, the evaluation of the total variation distance. By the triangle inequality, it suffices to show that $\|\mathbf{P}^n - \mathcal{T}^{n,1}\mathbf{P}^n\|_{\mathrm{TV}}$ and $\|\mathcal{T}^{n,2} \circ \mathcal{T}^{n,1}\mathbf{P}^n - \mathcal{T}^{n,1}\mathbf{P}^n\|_{\mathrm{TV}}$ can be made arbitrarily small for an appropriate choice of $I$ and for large enough $n$. Hence, we need to compare centered Gaussian measures. Let us start by evaluating the distance between the measures $\mathbf{P}^n$ and $\mathcal{T}^{n,1}\mathbf{P}^n$.

Recalling the construction of $\mathcal{X}^n$ in Section 7.2.1, these two measures only differ on the space of low frequencies $\mathcal{X}_\mathbf{e}^n$, and the covariance matrix of $\mathbf{P}^n$ on this space of dimension $m = 2^{j_0+2} + 1$ is the Toeplitz matrix $T_m(f_0)$ with the function $f_0(s) = |2\sin(\frac{s}{2})|^{1-2H_0}$ [the notation $T_m(f)$ is for the matrix with entries $T_m(f)_{k,l} := \frac{1}{2\pi} \int_{-\pi}^{\pi} f(s) e^{i(k-l)s} ds$ for $1 \leq k, l \leq m$]. The Gaussian measure $\mathcal{T}^{n,1}\mathbf{P}^n$ has, on the same space, covariance matrix $T_m(f_1)$ with $f_1(s) = |2\sin(\frac{s}{2})|^{1-2H_1}$. Then some considerations of Gaussian measures enable us to control, here, the distance between these two measures by the trace bound $\mathrm{Tr}([T_m(f_1)T_m(f_0)^{-1} - \mathrm{Id}]^2)$. Now the proof consists in making the following sequence of approximations rigorous:

$$\mathrm{Tr}([T_m(f_1)T_m(f_0)^{-1} - \mathrm{Id}]^2) \asymp \mathrm{Tr}\left(T_m\left[\left(\frac{f_1}{f_0} - 1\right)^2\right]\right)$$
$$\asymp \frac{m}{2\pi} \int_{-\pi}^{\pi} \left(\frac{f_1(s)}{f_0(s)} - 1\right)^2 ds$$
$$\asymp cm\varepsilon^2 \asymp c2^{j_0}\varepsilon^2 \asymp cI^{-2}.$$

The first approximation above expresses the quasi-homomorphism property of the Toeplitz operator $f \rightsquigarrow T_m(f)$, while the second one is a kind of Szegö theorem. The third approximation is obtained since $\|\frac{f_1}{f_0} - 1\|_2 \leq \varepsilon$, where the $L^2$-norm is taken over $[-\pi, \pi]$. Again, a detailed proof is presented in [11], where we use the method developed in Dahlhaus [5] and Fox and Taqqu [8] to deal with Toeplitz matrices (and Brockwell and Davis [3] too for more elementary results).

Finally, the control of $\|\mathcal{T}^{n,2} \circ \mathcal{T}^{n,1} \mathbf{P}^n - \mathcal{T}^{n,1} \mathbf{P}^n\|_{\mathrm{TV}}$ is obtained by similar techniques (see [11] for details). The property (iv) of Proposition 5 is proved.

**Acknowledgments.** We are grateful to Yuri Golubev, François Roueff and Philippe Soulier for helpful discussions and comments.

LABORATOIRE D'ANALYSE
ET DE MATHÉMATIQUE APPLIQUÉES
CNRS UMR 8050
UNIVERSITÉ DE MARNE-LA-VALLÉE
5 BOULEVARD DESCARTES
CHAMPS-SUR-MARNE 77454
FRANCE
E-MAIL: Arnaud.Gloter@univ-mlv.fr
        Marc.Hoffmann@univ-mlv.fr